\documentclass[default]{svmult}

\usepackage{type1cm}        
\usepackage{makeidx}         
\usepackage{graphicx}        
\usepackage{multicol}        
\usepackage[bottom]{footmisc}

\usepackage{newtxtext}       %
\usepackage{newtxmath}       

\makeindex             

\newcommand{\mumarg}[0]{\mu^{\text{M}}}

\begin{document}

\title*{Error bounds for some approximate posterior measures in Bayesian inference}
\author{Han Cheng Lie and T. J. Sullivan and Aretha Teckentrup}
\institute{Han Cheng Lie \at Institut f\"{u}r Mathematik, Universit\"{a}t Potsdam, Campus Golm, Haus 9, Karl-Liebknecht Strasse 24-25, Potsdam OT Golm 14476, Germany, \email{hanlie@uni-potsdam.de}
\and T. J. Sullivan  \at  Mathematics Institute and School of Engineering, The University of Warwick, Coventry CV4 7AL, United Kingdom, \email{t.j.sullivan@warwick.ac.uk}, and Zuse-Institut Berlin, Takustrasse 7, Berlin 14195, Germany, \email{sullivan@zib.de}
\and
Aretha Teckentrup \at School of Mathematics, University of Edinburgh, James Clerk Maxwell Building, Edinburgh, EH9 3FD, United Kingdom, \email{a.teckentrup@ed.ac.uk}}
%
%
\maketitle

\abstract*{In certain applications involving the solution of a Bayesian inverse problem, it may not be possible or desirable to evaluate the full posterior, e.g. due to the high computational cost of doing so. This problem motivates the use of approximate posteriors that arise from approximating the data misfit or forward model. We review some error bounds for random and deterministic approximate posteriors that arise when the approximate data misfits and approximate forward models are random.}

\abstract{In certain applications involving the solution of a Bayesian inverse problem, it may not be possible or desirable to evaluate the full posterior, e.g. due to the high computational cost of doing so. This problem motivates the use of approximate posteriors that arise from approximating the data misfit or forward model. We review some error bounds for random and deterministic approximate posteriors that arise when the approximate data misfits and approximate forward models are random.}

\section{Introduction}
\label{sec:intro}
An inverse problem consists of recovering an unknown parameter $u$ that belongs to a possibly infinite-dimensional space $\mathcal{U}$ from noisy data $y$ of the form
\begin{equation}
\label{eq:observation_model}
y=G(u)+\eta\in\mathcal{Y},
\end{equation}
where $\mathcal{Y}$ is the `data space', $G:\mathcal{U}\to\mathcal{Y}$ is a known `forward operator', and $\eta$ is a random variable. In many problems of interest, the parameter space $\mathcal{U}$ is a subset of an infinite-dimensional Banach space, the data space $\mathcal{Y}$ is often taken to be $\mathbb{R}^d$ for some possibly large $d\in\mathbb{N}$, and $\eta$ is assumed to be Gaussian.

One of the main difficulties with inverse problems is that they often do not satisfy Hadamard's definition of well-posedness. To circumvent this difficulty, one may use the Bayesian approach, in which one incorporates information about the unknown $u$ from existing data and from new data in the `prior' probability measure $\mu_0$ on $\mathcal{U}$ and in the `data misfit' $\Phi:\mathcal{Y}\times\mathcal{U}\to\mathbb{R}$ respectively. If $\eta\in\mathbb{R}^d$ in \eqref{eq:observation_model} is distributed according to the normal distribution $N(0,\Gamma)$ with positive definite $\Gamma\in\mathbb{R}^{d\times d}$, then 
\begin{equation}
\label{eq:quadratic_potential}
\Phi(y,u):=\frac{1}{2}\Vert \Gamma^{-1/2}\left(y-G(u)\right)\Vert^2.
\end{equation}
By Bayes' formula, the posterior $\mu^y$ is a probability measure on $\mathcal{U}$ that is absolutely continuous with respect to the prior $\mu_0$, and has Radon--Nikodym derivative 
\begin{equation}
\label{eq:true_posterior}
\frac{\mathrm{d}\mu^y}{\mathrm{d}\mu_0}(u):=\frac{\exp(-\Phi(y,u))}{Z(y)},\quad Z(y):=\int_{\mathcal{U}}\exp(-\Phi(y,u'))\mathrm{d}\mu_0(u').
\end{equation}
The posterior $\mu^y$ describes the distribution of the unknown $u$, conditioned upon the data $y$. By imposing conditions jointly upon $\Phi$ and $\mu_0$, one can show that the Bayesian solution $\mu^y$ to the inverse problem depends continuously on the data, and one can prove the well-posedness of the Bayesian inverse problem; see \cite{DashtiStuart2017}.

For simplicity, we shall assume that the data $y$ is given and fixed, and omit the dependence of the posterior, data misfit, and normalisation constant $Z$ on $y$. 

One challenge with solving Bayesian inverse problems in practice is that it is often not possible or desirable to evaluate the data misfit $\Phi(u)$ exactly. It then becomes necessary to find approximations $\Phi_N$ of the true data misfit $\Phi$ that can be computed more efficiently, such that for sufficiently large values of $N$, inference using the approximate misfit $\Phi_N$ effectively approximates inference using the true misfit $\Phi$. Thus, one needs to identify conditions on $\Phi_N$ such that two criteria are fulfilled: first, an approximate posterior measure $\mu_N$ defined by 
\begin{equation}
\label{eq:approximate_posterior}
\frac{\mathrm{d}\mu_N}{\mathrm{d}\mu_0}(u):=\frac{\exp(-\Phi_N(u))}{Z_N},\quad Z_N:=\int_{\mathcal{U}}\exp(-\Phi_N(u'))\mathrm{d}\mu_0(u')
\end{equation}
exists and is well-defined; and second, the approximate posterior $\mu_N$ provides an increasingly good approximation of the true posterior $\mu$ as the approximation parameter $N$ increases. In this paper, we review results from \cite{LST2018} that guarantee well-definedness of $\mu_N$ and establish error bounds for $\mu_N$ in terms of error bounds for $\Phi_N$.

In recent years, randomised numerical methods have been developed in order to overcome limitations of their deterministic counterparts. The field of probabilistic numerical methods \cite{Cockayne2019} injects randomness into existing deterministic solvers for differential equations in order to model the uncertainty due to unresolved subgrid-scale dynamics. Random approximations of the forward model have been applied for forward uncertainty propagation in a range of applications; see e.g. \cite{MarzoukXiu2009,Birolleau2014}.

Randomisation has been shown to yield gains in computational efficiency. Results from \cite[Section 5.7]{CF2005} showed a reduction by a factor of almost 25 in the CPU time needed for generating an independent sample with the Metropolis-Hastings algorithm, while in \cite{DKST2019}, a multilevel Markov Chain Monte Carlo method uses randomisation in the form of control variates for variance reduction. Stochastic programming ideas were used for more efficient posterior sampling in \cite{WBG2018}. The results we describe provide theoretical support for the use of randomisation in Bayesian inference, and extend the pioneering results from \cite{StuartTeckentrup2017}, which concerned Gaussian process approximations of data misfits and forward models.

To motivate the use of random approximate misfits, consider the following example: Let $X$ be any $\mathbb{R}^d$-valued random variable such that $\mathbb{E}[X]=0$ and $\mathbb{E}[XX^\top]$ is the $d\times d$ identity matrix, and let $\{X_i\}_{i\in\mathbb{N}}$ be i.i.d. copies of $X$. Given \eqref{eq:quadratic_potential}, 
\begin{align*}
 \Phi(u)=&\frac{1}{2}\left(\Gamma^{-1/2}(y-G(u)\right)\mathbb{E}\left[XX^\top\right]\left(\Gamma^{-1/2}(y-G(u)\right)
 \\
 =&\frac{1}{2}\mathbb{E}\left[\left\vert X^\top\left(\Gamma^{-1/2}(y-G(u))\right)\right\vert\right]\approx\frac{1}{2N}\sum_{j=1}^{N}\left\vert X_j^\top\left(\Gamma^{-1/2}(y-G(u))\right)\right\vert=:\Phi_N(u).
\end{align*}
In \cite{Le2017}, the misfit $\Phi_N$ above was used to obtain computational cost savings when solving inverse problems associated to PDE boundary value problems. The results we present below can be specialised to the case of $X$ with bounded support \cite[Proposition 4.1]{LST2018}. For example, we can use the \emph{$\ell$-sparse distribution} for some $0\leq \ell<1$; for $\ell=0$, this is the Rademacher distribution. Similar ideas have been applied for full waveform inversion in seismic tomography \cite{Aravkin2012}, for example.


\section{Error bounds for approximate posteriors}
\label{sec:error_bounds}

In what follows, we shall assume that the parameter space $\mathcal{U}$ admits a Borel $\sigma$-algebra, and we shall denote by $\mathcal{M}_1(\mathcal{U})$ the set of Borel probability measures on $\mathcal{U}$. Recall that the Hellinger metric $d_{\mathrm{H}}:\mathcal{M}_1(\mathcal{U})\times\mathcal{M}_1(\mathcal{U})\to [0,1]$ is defined by
\begin{equation*}
d_{\mathrm{H}}(\mu,\nu)^2:=\frac{1}{2}\int_{\mathcal{U}}\left\vert \sqrt{\frac{\mathrm{d}\mu}{\mathrm{d}\pi}}(u')-\sqrt{\frac{\mathrm{d}\nu}{\mathrm{d}\pi}}(u')\right\vert^2\mathrm{d}\pi(u'),
\end{equation*}
where $\pi\in\mathcal{M}_1(\mathcal{U})$ is any measure such that $\mu$ and $\nu$ are both absolutely continuous with respect to $\pi$. It is known that $d_{\mathrm{H}}$ does not depend on the choice of $\pi$.

\subsection{Error bounds for random approximate posteriors}

We first present error bounds on random approximate posteriors $\mu_N$ associated to random misfits $\Phi_N$, where $N\in\mathbb{N}$. That is, given a probability space $(\Omega,\mathcal{F},\mathbb{P})$, we shall view a random misfit as a measurable function $\Phi_N:\Omega\times\mathcal{U}\to\mathbb{R}$. Furthermore, we shall assume that the randomness associated to the approximate misfit $\Phi_N$ is independent of the randomness associated to the unknown parameter $u$. In what follows, $\nu_N$ denotes a probability measure on $\Omega$ with the property that the distribution of the random function $\Phi_N$ is given by $\nu_N\otimes \mu_0$. 

Given \eqref{eq:true_posterior} and \eqref{eq:approximate_posterior}, a natural question is to establish an appropriate bound on the Hellinger distance between the true posterior $\mu$ and the approximate posterior $\mu_N$ in terms of some norm of the error between the true misfit $\Phi$ and the approximate misfit $\Phi_N$. We emphasise that the approximate posterior $\mu_N$ in \eqref{eq:approximate_posterior} is random in the sense that it depends on $\omega$, since the approximate misfit $\Phi_N$ depends on $\omega$. Therefore, the Hellinger distance $d_{\mathrm{H}}(\mu,\mu_N)$ will depend on $\omega$ as well. To describe such a bound, we shall take the expectation of the Hellinger distance with respect to $\nu_N$, and let
\begin{equation*}
\left\Vert \mathbb{E}_{\nu_N}\left[f\left(\Phi_N\right)\right]~\right\Vert_{L^{q}_{\mu_0}(\mathcal{U})}:=\left(\int_{\mathcal{U}} \left\vert\int_{\Omega} f\left(\Phi_N(\omega,u)\right)\mathrm{d}\nu_N(\omega)\right\vert^q\mathrm{d}\mu_0(u)\right)^{1/q}
\end{equation*}
for any Borel-measurable function $f:\mathbb{R}\to\mathbb{R}$ and $q\in [1,\infty)$. We define the quantity $\Vert \mathbb{E}_{\nu_N}[f(\Phi_N)]\Vert_{L^\infty_{\mu_0}(\mathcal{U})}$ analogously. With these preparations, we present the following theorem, which was given in \cite[Theorem 3.2]{LST2018}.
\begin{theorem}[Error bound for random approximate posterior]
\label{thm:1}
Let $(q_1,q_1')$ and $(q_2,q_2')$ be pairs of H\"{o}lder conjugate exponents, and let $D_1$, $D_2$ be positive scalars that depend only on $q_1$ and $q_2$. Suppose the following conditions hold:
\begin{eqnarray}
\left\Vert\mathbb{E}_{\nu_N}\left[\left(\exp\left(-\tfrac{1}{2}\Phi\right)+\exp\left(-\tfrac{1}{2}\Phi_N\right)\right)^{2q_1}\right]^{1/q_1}\right\Vert_{L^{q_2}_{\mu_0}(\mathcal{U})}\leq& D_1
\label{eq:A1}
\\
\left\Vert \mathbb{E}_{\nu_N}\left[\left(Z_N\max\{Z^{-3},Z_N^{-3}\}\left(\exp\left(-\Phi\right)+\exp\left(-\Phi_N\right)\right)^{2}\right)^{q_1}\right]^{1/q_1}\right\Vert_{L^{q_2}_{\mu_0}(\mathcal{U})}\leq& D_2.
\label{eq:A2}
\end{eqnarray}
Then 
\begin{equation*}
\mathbb{E}_{\nu_N}\left[d_{\mathrm{H}}\left(\mu,\mu_N\right)^2\right]^{1/2}\leq (D_1+D_2)\left\Vert \mathbb{E}_{\nu_N}\left[\left\vert \Phi-\Phi_N\right\vert^{2q'_1}\right]^{1/2q_1'}\right\Vert_{L^{2q_2'}_{\mu_0}(\mathcal{U})}.
\end{equation*}
\end{theorem}
Theorem \ref{thm:1} provides a bound on the mean square Hellinger distance between the true posterior $\mu$ and the random approximate posterior $\mu_N$, in terms of an appropriate norm of the error $\Phi-\Phi_N$. The bound \eqref{eq:A1} implies that the negative tails of both $\Phi$ and $\Phi_N$ must decay exponentially quickly with respect to the $\nu_N\otimes \mu_0$-measure, and is satisfied, for example, when both $\Phi$ and $\Phi_N$ are bounded from below. Since $Z_N\max\{Z^{-3},Z_N^{-3}\}=\max\{Z_NZ^{-3},Z_N^{-2}\}$, it follows that the constraint imposed on the misfit $\Phi_N$ by \eqref{eq:A2} is that $\exp(-\Phi_N)$ should be neither too concentrated nor too broad. Together, conditions \eqref{eq:A1} and \eqref{eq:A2} ensure that the random approximate posterior $\mu_N$ exists, is well-defined, and satisfies the desired bound on the mean square Hellinger distance with respect to the true posterior $\mu$.

An alternative way to generate an approximate posterior measure given a random approximate misfit is to compute a marginal approximate posterior $\mumarg_N$, defined by
\begin{equation}
\label{eq:marginal_approximate_posterior}
\frac{\mathrm{d}\mumarg_N}{\mathrm{d}\mu_0}(u):=\frac{\mathbb{E}_{\nu_N}\left[\exp(-\phi_N(u))\right]}{\mathbb{E}_{\nu_N}\left[Z_N\right]}.
\end{equation}
Note that, since we have taken expectations with respect to $\nu_N$, the marginal approximate posterior does not depend on $\omega$, and is in this sense deterministic. The following theorem was given in \cite[Theorem 3.1]{LST2018}.
\begin{theorem}[Error bound for marginal approximate posterior]
\label{thm:2}
Let $(p_1,p_1')$, $(p_2,p_2')$, and $(p_3,p_3')$ be H\"{o}lder conjugate exponent pairs, and suppose there exist finite, positive scalars $C_1$, $C_2$, and $C_3$ that depend only on $p_1$, $p_2$, and $p_3$, such that the following conditions hold:
\begin{eqnarray}
\min\left\{\left\Vert \mathbb{E}_{\nu_N}\left[\exp\left(-\Phi_N\right)\right]^{-1}\right\Vert_{L^{p_1}_{\mu_0}(\mathcal{U})},\left\Vert\exp(\Phi)\right\Vert_{L^{p_1}_{\mu_0}(\mathcal{U})}\right\}\leq& C_1
\label{eq:B1}
\\
\left\Vert \mathbb{E}_{\nu_N}\left[\left(\exp(-\Phi)+\exp(-\Phi_N)\right)^{p_2}\right]^{1/p_2}\right\Vert_{L^{2p_1'p_3}_{\mu_0}(\mathcal{U})}\leq& C_2
\label{eq:B2}
\\
C_3^{-1}\leq \mathbb{E}_{\nu_N}\left[Z_N\right]\leq C_3.
\label{eq:B3}
\end{eqnarray}
Then there exists $C>0$ that does not depend on $N$ such that
\begin{equation*}
d_{\mathrm{H}}(\mu,\mumarg_N)\leq C\left\Vert \mathbb{E}_{\nu_N}\left[\left\vert \Phi-\Phi_N\right\vert^{p_2'}\right]^{1/p_2'}\right\Vert_{L^{2 p_1'p_3'}_{\mu_0}(\mathcal{U})}.
\end{equation*}
\end{theorem}
The bounds in \eqref{eq:B3} ensure that the denominator in \eqref{eq:marginal_approximate_posterior} is strictly positive and finite. Thus, these bounds play a fundamental role in ensuring that the marginal approximate posterior exists and is well-defined. The bound in \eqref{eq:B2} reiterates the bound \eqref{eq:A1}, modulo the $\tfrac{1}{2}$ factor, and thus serves a similar purpose as \eqref{eq:A1}. The bound in \eqref{eq:B1} serves a similar purpose as \eqref{eq:A2}. However, the minimum operator implies that it is not necessary for both $\Phi$ and $\Phi_N$ to be well-behaved. 

The following result is a corollary of Theorem \ref{thm:1}, Theorem \ref{thm:2}, and \cite[Lemma 3.5]{LST2018}. The main idea is to specify sufficient conditions for the hypotheses of both Theorem \ref{thm:1} and Theorem \ref{thm:2} to hold.
\begin{corollary}[Joint conditions for error bounds on both approximate posteriors]
\label{cor:01}
Suppose the following conditions are satisfied:
\begin{enumerate}
\item[(i)] There exists $C_0\in\mathbb{R}$ that does not depend on $N$ such that $\Phi\geq -C_0$ on $\mathcal{U}$ and, for all $N\in\mathbb{N}$, $\nu_N(\Phi_N\geq -C_0)=1$,
\item[(ii)] For any $0<C_3<\infty$ such that $C_3^{-1}<Z<C_3$, there exists $N^\ast(C_3)\in\mathbb{N}$ such that $N\geq N^\ast$ implies 
\begin{equation*}
\left\Vert \mathbb{E}_{\nu_N}[\vert \Phi-\Phi_N\vert]\right\Vert_{L^1_{\mu_0}(\mathcal{U})}\leq \frac{1}{2}\exp(-C_0)\min\left\{Z-C_3^{-1},C_3-Z\right\},
\end{equation*}
and
\item[(iii)] there exists some $2<\rho^\ast<+\infty$ such that $\Vert \mathbb{E}_{\nu_N}[\exp(\rho^\ast\Phi_N)]\Vert_{L^1_{\mu_0}(\mathcal{U})}$ is finite.
\end{enumerate}
Then for each $N\geq N^\ast(C_3)$,
\begin{equation}
\label{eq:corollary_conclusion1}
d_{\mathrm{H}}\left(\mu,\mumarg_N\right)\leq C\left\Vert \mathbb{E}_{\nu_N}\left[\left\vert\Phi-\Phi_N\right\vert\right]\right\Vert_{L^{2\rho^\ast/(\rho^\ast-1)}_{\mu_0}(\mathcal{U})}
\end{equation}
and 
\begin{equation}
\label{eq:corollary_conclusion2}
\mathbb{E}_{\nu_N}\left[d_{\mathrm{H}}(\mu,\mu_N)^2\right]^{1/2}\leq D\left\Vert \mathbb{E}_{\nu_N}\left[\left\vert \Phi-\Phi_N\right\vert^{2\rho^\ast/(\rho^\ast-2)}\right]^{(\rho^\ast-2)/(2\rho^\ast)}\right\Vert_{L^{1}_{\mu_0}(\mathcal{U})},
\end{equation}
where $C,D>0$ depend on $\Vert \mathbb{E}_{\nu_N}[\exp(\rho^\ast\Phi_N)]\Vert_{L^1_{\mu_0}(\mathcal{U})}^{1/\rho^\ast}$. If in addition to conditions $(i)$--$(iii)$ it holds that 
\begin{equation*}
\sup_{N\geq N^\ast(C_3)}\left\Vert \mathbb{E}_{\nu_N}\left[\exp(\rho^\ast\Phi_N)\right]\right\Vert_{L^1_{\mu_0}(\mathcal{U})}<\infty,
\end{equation*}
then the constants $C$ and $D$ in \eqref{eq:corollary_conclusion1} and \eqref{eq:corollary_conclusion2} do not depend on $N$.
\end{corollary}

Condition $(i)$ amounts to a common uniform lower bound on all the misfits, both the true misfit and the collection of random approximate misfits, and thus plays a role in ensuring that \eqref{eq:A1} and \eqref{eq:B2} are satisfied. Condition $(ii)$ makes precise the assumption that $\Phi_N$ approximates $\Phi$ in the $L^1_{\nu_N\otimes \mu_0}$ topology, which is a necessary condition for ensuring that the right-hand sides of the conclusions of Theorem \ref{thm:1} and Theorem \ref{thm:2} are finite. Condition $(iii)$ describes an exponential integrability condition on the random approximate misfits and ensures that \eqref{eq:A2} and \eqref{eq:B1} are satisfied. Thus the additional condition amounts to a uniform exponential integrability condition over all sufficiently large values of $N$.

\begin{remark}
\label{rem:boundedness_from_below_of_negative_log_likelihoods}
Neither Theorem \ref{thm:1} nor Theorem \ref{thm:2} require boundedness from below of either $\Phi$ or the $\Phi_N$. However, the negative tails of both $\Phi$ and $\Phi_N$ must decay exponentially quickly at a sufficiently high rate, as specified by \eqref{eq:B2} and \eqref{eq:A1} respectively.
\end{remark}

\subsection{Error bounds for random forward models}

Next, we consider approximate posterior measures that arise as a result of approximating the forward model $G$ in \eqref{eq:observation_model}. For simplicity, we shall consider only the case when the data misfit $\Phi$ and forward model $G$ are related via the quadratic potential \eqref{eq:quadratic_potential}. In particular, this means that if $G_N:\mathcal{U}\to\mathcal{Y}$ is an approximation of the true forward model $G$, then the resulting approximate data misfit is given by 
\begin{equation*}
\Phi_N(u):=\frac{1}{2}\Vert \Gamma^{-1}\left(y-G_N(u)\right)\Vert^2.
\end{equation*}
The following theorem is a nonasymptotic reformulation of \cite[Theorem 3.9 (b)]{LST2018}. 
\begin{theorem}[Error bounds for approximate posteriors] 
Suppose there exists $2<\rho^\ast<\infty$ such that $\sup_N\mathbb{E}_{\nu_N}[\exp(\rho^\ast\Phi_N)]\in L^1_{\mu_0}(\mathcal{U})$ is finite. If there exists an $N^\ast\in\mathbb{N}$ such that, for all $N\geq N^\ast$,
\begin{equation*}
\left\Vert \mathbb{E}_{\nu_N}\left[\left\Vert G-G_N\right\Vert^{4\rho^\ast/(\rho^\ast-2)}\right]^{(\rho^\ast-2)/(2\rho^\ast)}\right\Vert_{L^{2\rho^\ast/(\rho^\ast-1)}_{\mu_0}(\mathcal{U})}\leq 1,
\end{equation*}
then 
\begin{equation*}
d_{\mathrm{H}}\left(\mu,\mumarg_N\right)\leq C\left\Vert \mathbb{E}_{\nu_N}\left[\left\Vert G_N-G\right\Vert^{2}\right]\right\Vert^{1/2}_{L^{2\rho^\ast/(\rho^\ast-1)}_{\mu_0}(\mathcal{U})}
\end{equation*}
and 
\begin{equation*}
\mathbb{E}_{\nu_N}\left[d_{\mathrm{H}}(\mu,\mu_N))^2\right]^{1/2}\leq D\left\Vert \mathbb{E}_{\nu_N}\left[\left\Vert G_N-G\right\Vert^{4\rho^\ast/(\rho^\ast-2)}\right]^{(\rho^\ast-2)/(2\rho^\ast)}\right\Vert^{1/2}_{L^{2}_{\mu_0}(\mathcal{U})}
\end{equation*}
for $C,D>0$ that do not depend on $N$.
\end{theorem}
The theorem can be rewritten so that, instead of imposing a uniform exponential integrability condition on the approximate quadratic potentials $\Phi_N$, one instead imposes an exponential integrability condition on the true data misfit $\Phi$; see \cite[Theorem 3.9 (a)]{LST2018}. An additional hypothesis in this case is that the expectations of the approximate data misfit functions are $\nu_N$-almost surely bounded, in the sense that $\nu_N(\Phi_N\ \vert \ \mathbb{E}_{\mu_0}[\Phi_N]\leq C_4)=1$ for some $C_4\in\mathbb{R}$ that does not depend on $N$. 

\section{Conclusions and directions for future work}
\label{sec:conclusion}

This paper has reviewed the main error bounds of \cite{LST2018} concerning deterministic and random approximate posteriors that arise when performing Bayesian inference with random approximate data misfits or random forward models. The error bounds on the approximate posterior measures are given with respect to the Hellinger metric on the space of Borel probability measures $\mathcal{M}_1(\mathcal{U})$. Given a fixed prior measure $\mu_0$, these error bounds describe -- with specific exponents of integrability and problem-dependent constants -- the local or global Lipschitz continuity of the map that takes a data misfit as input and produces the corresponding posterior measure as output. Aside from the regularity assumptions made on the random approximations, the error bounds shown above make no structural assumptions on the approximations used. For example, we do not assume that the random approximations involve Gaussian random variables, or random variables with bounded support.

Recent work has highlighted the importance of considering other metrics on $\mathcal{M}_1(\mathcal{U})$, and also of proving well-posedness of the solution of a Bayesian inverse problem by establishing continuous (instead of Lipschitz continuous) dependence on either the data, prior, or data misfit. The well-posedness of Bayesian inverse problems in the sense of continuous dependence with respect to the data of the posterior for given prior and data misfit was established in \cite{Latz2019}. Local Lipschitz continuity with respect to \emph{deterministic} perturbations in the prior or data misfit was shown in \cite{Sprungk2019}. In both \cite{Latz2019,Sprungk2019}, continuity is with respect to the topologies induced by the total variation metric, by Wasserstein $p$-metrics, or by the Kullback-Leibler divergence. 

A key assumption made in \cite{Sprungk2019} when establishing local Lipschitz continuity for a fixed prior $\mu_0$ with respect to perturbations in the data misfit is that the deterministic perturbed data misfit is $\mu_0$-almost surely bounded from below. As highlighted in Remark \ref{rem:boundedness_from_below_of_negative_log_likelihoods}, the analysis of \cite{LST2018} does not require that either the true data misfit or the random approximate log-likelihood are $\mu_0$-almost surely bounded from below. For future work, we aim to establish similar continuity results with respect to different metrics, as demonstrated in \cite{Latz2019,Sprungk2019}, but at the same level of generality of \cite{LST2018}.

\begin{acknowledgement}
The research of HCL has been partially funded by Deutsche Forschungsgemeinschaft (DFG) - SFB1294/1 - 318763901 and by Universit\"{a}t Potsdam. The work of TJS has been partially supported by the Freie Universit\"{a}t Berlin within the Excellence Initiative of the German Research Foundation (DFG). The authors thank an anonymous referee for their feedback.
\end{acknowledgement}
\biblstarthook{
}

\end{document}